\documentclass[twoside,10pt]{article}
\usepackage{graphicx}
\usepackage{subfigure}
\usepackage{url}
\usepackage{color}

\usepackage{amsfonts}
\newcommand\R{\mathbb R}
\newcommand\T{\mathbb T}
\newcommand\Z{\mathbb Z}
\newcommand\ee{\mathrm e}

\newcommand\eps\varepsilon
\newcommand\A{\mathcal A}
\newcommand\J{\mathcal J}
\newcommand\Lc{\mathcal L}
\newcommand\M{\mathcal M}
\newcommand\Rc{\mathcal R}
\newcommand\W{\mathcal W}
\newcommand\const{{\rm const}}

\newcommand\ds{\displaystyle}
\newcommand\dfrac[2]{\ds\frac{#1}{#2}}

\newcommand\p[1]{\left(#1\right)}
\newcommand\pq[1]{\left[#1\right]}
\newcommand\pp[1]{\left\{#1\right\}}
\newcommand\scprod[2]{\left\langle#1,#2\right\rangle}
\newcommand\abs[1]{\left|#1\right|}

\newcommand\tl{\tilde}
\newcommand\wtl{\widetilde}
\newcommand\wh[1]{\widehat{#1}}
\newcommand\ol[1]{\overline{#1}}
\newcommand\ool[1]{\ol{\ol{#1}}}
\newcommand\df{{\rm d}}
\newcommand\Df{{\rm D}}
\newcommand\Ord{\mathcal O}
\newcommand\rint{\mathop{\rm rint}}
\newcommand\ut{{\rm u}}  
\newcommand\st{{\rm s}}  
\newcommand\hardsum[2]
  {\ds\sum_{\begin{array}{c}\null^{#1}\\[-3pt]\null^{#2}\end{array}}}
\newcommand\tr{\mathop{\rm tr}}

\newtheorem{theorem}{Theorem}

\newtheorem{lemma}[theorem]{Lemma}
\newcommand\bremark{\noindent{\bf Remark.}\ \ }
\newcommand\eremark{\bigskip}
\newcommand\proof{\noindent\emph{Proof.}\ \ }
\newcommand\proofof[1]{\noindent\emph{Proof of #1.}\ \ }
\newcommand\qed{\ \ \null\nolinebreak\hfill$\frame{\large\phantom a}$}

\newcommand\beq{\begin{equation}}
\newcommand\eeq{\end{equation}}
\newcommand\bea{\begin{eqnarray}}
\newcommand\eea{\end{eqnarray}}
\newcommand\bean{\begin{eqnarray*}}
\newcommand\eean{\end{eqnarray*}}
\newcommand\btm{\vspace{-\baselineskip}\begin{itemize}}
\newcommand\etm{\end{itemize}\vspace{-\baselineskip}}

\begin{document}

\title{Continuation of the exponentially small transversality for the splitting
  of separatrices to a whiskered torus with silver ratio
  \ \footnote{This work has been partially supported by the Spanish
      MINECO-FEDER Grant MTM2012-31714,
      the Catalan Grant 2014SGR504,
      and the Russian Scientific Foundation Grant 14-41-00044.
      The author~MG has also been supported by the
      DFG~Collaborative Research Center TRR~109
      ``Discretization in Geometry and Dynamics''.}}
\author{\sc
    Amadeu Delshams$\,^1$,
  \ Marina Gonchenko$\,^2$,
\\[4pt]\sc
  Pere Guti\'errez$\,^1$
\\[12pt]
  {\small
  $^1\;$\parbox[t]{5.2cm}{
    Dep. de Matem\`atica Aplicada I\\
    Universitat Polit\`ecnica de Catalunya\\
    Av. Diagonal 647, 08028 Barcelona\\
    {\footnotesize
      \texttt{amadeu.delshams@upc.edu}\\
      \texttt{pere.gutierrez@upc.edu}}}
  \quad
  $^2\;$\parbox[t]{5.2cm}{
    Technische Universit\"at Berlin\\
    Institut f\"ur Mathematik\\
    Stra{\ss}e des 17. Juni 136\\
    D-10623 Berlin\\
    {\footnotesize
      \texttt{gonchenk@math.tu-berlin.de}}}
  }}
\maketitle
\begin{abstract}
We study the exponentially small splitting of invariant manifolds of whiskered
(hyperbolic) tori with two fast frequencies in nearly-integrable Hamiltonian
systems whose hyperbolic part is given by a pendulum. We consider a torus whose
frequency ratio is the silver number $\Omega=\sqrt{2}-1$. We show that the
Poincar\'e--Melnikov method can be applied to establish the existence of
4 transverse homoclinic orbits to the whiskered torus, and provide asymptotic
estimates for the tranversality of the splitting whose dependence on the
perturbation parameter $\varepsilon$ satisfies a periodicity property. We also
prove the continuation of the transversality of the homoclinic orbits for all
the sufficiently small values of~$\varepsilon$, generalizing the results
previously known for the golden number.
\par\vspace{12pt}
\end{abstract}

\section{Introduction and setup}

\subsection{Background and state of the art}

This paper is dedicated to the study of the transversality of the exponentially
small splitting of separatrices in a perturbed
3-degree-of-freedom Hamiltonian system,
associated to a 2-dimensional whiskered torus
(invariant hyperbolic torus) whose frequency ratio
is the silver number $\Omega=\sqrt{2}-1$.
This quadratic irrational number has nice arithmetic properties since it has a
1-periodic continued fraction.

We start with an integrable Hamiltonian $H_0$
having whiskered (hyperbolic) tori with a \emph{separatrix}:
coincident stable and unstable whiskers (invariant manifolds).
We focus our attention on a torus,
with a frequency vector of \emph{fast frequencies}:
\beq
  \omega_\eps = \frac\omega{\sqrt{\eps}}\;,
  \qquad
  \omega=(1,\Omega), \qquad \Omega=\sqrt{2}-1.
\label{eq:omega_eps}
\eeq
This frequency ratio $\Omega$ is called the \emph{silver number}.
If we consider a perturbed Hamiltonian $H=H_0+\mu H_1$, where $\mu$ is small,
in general the stable and unstable whiskers do not coincide anymore,
and this phenomenon has got the name of \emph{splitting of separatrices}.
If we assume, for the two involved parameters,
a relation of the form $\mu=\eps^p$
for some $p>0$, we have a problem of singular perturbation
and in this case the splitting is \emph{exponentially small}
with respect to $\eps$.
Our aim is to detect homoclinic orbits associated
to persistent whiskered tori,
provide \emph{asymptotic estimates} for both the splitting distance
and its \emph{transversality}, and
{use the arithmetic properties of the silver number $\Omega$ in order to show
the \emph{continuation} of the transversality of
the homoclinic orbits for \emph{all} sufficiently small $\eps$}.
When transversality takes place, the perturbed system turns out to
be non-integrable and there is chaotic dynamics near the homoclinic orbits.

A very usual tool to measure the splitting is
the \emph{Poincar\'e--Melnikov method}, introduced by
Poincar\'e in \cite{Poincare90} and rediscovered much later by
Melnikov and Arnold \cite{Melnikov63,Arnold64}.
By considering a transverse section to the
stable and unstable perturbed whiskers,
one can consider a function $\M(\theta)$, $\theta\in\T^2$,
usually called the \emph{splitting function},
giving the vector distance between the whiskers on this section.
The method provides a first order approximation
to this function, with respect to the parameter $\mu$,
given by the \emph{Melnikov function} $M(\theta)$, defined by an integral.
We have
\beq\label{eq:melniapproxM}
  \M(\theta)=\mu M(\theta)+\Ord(\mu^2),
\eeq
and hence for $\mu$ small enough the simple zeros $\theta_*$ of $M(\theta)$
give rise to transverse intersections between the perturbed whiskers.
In this way, we can obtain asymptotic estimates for both the
\emph{maximal splitting distance} as the maximum
of the function $\abs{\M(\theta)}$,
and for the \emph{transversality} of the splitting,
which can be measured by the minimal eigenvalue (in modulus)
of the ($2\times2$)-matrix $\Df\M(\theta_*)$.

An important related fact is that both functions $\M(\theta)$ and $M(\theta)$
are gradients of scalar functions \cite{Eliasson94,DelshamsG00}:
\[
  \M(\theta)=\nabla\Lc(\theta),
  \qquad
  M(\theta)=\nabla L(\theta).
\]
Such scalar funtions are called \emph{splitting potential} and
\emph{Melnikov potential} respectively, and
the transverse homoclinic orbits correspond
to the nondegenerate critical points of the splitting potential.

As said before, the case of fast frequencies $\omega_\eps$
as in~(\ref{eq:omega_eps}), with a perturbation of order $\mu=\eps^p$,
turns out to be a \emph{singular problem}.
The difficulty is that the Melnikov function $M(\theta)$ is exponentially small
in $\eps$, and the Poincar\'e--Melnikov method cannot be directly applied,
unless we assume that $\mu$ is exponentially small with respect to $\eps$.
In order to validate the method in the case $\mu=\eps^p$,
with $p$ as small as possible, it was introduced in \cite{Lazutkin03}
the use of parameterizations of a complex strip of the whiskers
(whose width is defined by the singularities of the unperturbed ones),
together with flow-box coordinates,
in order to ensure that the error term is also exponentially small,
and that the Poincar\'e-Melnikov approximation dominates it.
This tool was initially developed for the Chirikov standard map
\cite{Lazutkin03}, for Hamiltonians with one and a half degrees of freedom
(with 1~frequency) \cite{DelshamsS92,DelshamsS97,Gelfreich97}
and for area-preserving maps \cite{DelshamsR98}.

Later, those methods were extended to the case of whiskered tori
with 2~frequencies. In this case, the arithmetic properties of the frequencies
play an important role in the exponentially small asymptotic estimates of the
splitting function, due to the presence of \emph{small divisors}.
This was first mentioned in \cite{Lochak90}
and later detected in \cite{Simo94}, and then rigorously proved in
\cite{DelshamsGJS97} for the quasi-periodically forced pendulum,
assuming a polynomial perturbation in the coordinates associated
to the pendulum. Recently, a more general (meromorphic) perturbation
has been considered in \cite{GuardiaS12}.
It is worth mentioning that, in some cases, the Poincar\'e--Melnikov method
does not predict correctly
the size of the splitting, as shown in \cite{BaldomaFGS12}.

As an alternative way to study the splitting, the parametrization of the
whiskers as solutions of the Hamilton--Jacobi equation was used in
\cite{Sauzin01,LochakMS03,RudnevW00},
and exponentially small estimates are also obtained with this method,
as well as the transversality of the splitting,
provided some intervals of the perturbation parameter~$\eps$ are excluded.
Similar results were obtained in \cite{DelshamsG04,DelshamsG03}.
Besides, in the case of golden ratio $\Omega=(\sqrt5-1)/2$ it was
shown in \cite{DelshamsG04} the \emph{continuation}
of the transversality for \emph{all} sufficiently small values of $\eps$,
under a certain condition on the phases of the perturbation.
Otherwise, homoclinic bifurcations can occur,
studied, {for instance}, in \cite{SimoV01} for the Arnold's example.
The generalization of this approach to some other quadratic frequency ratios
was considered in \cite{DelshamsG03},
extending the asymptotic estimates for the splitting,
but without a satisfactory result concerning
the continuation of the transversality.
Recently, a parallel study for  the cases of~2 and 3~frequencies has been
considered in \cite{DelshamsGG14a}
(in the case of 3~frequencies, with a frequency vector
$\omega=(1,\Omega,\Omega^2)$, where $\Omega$ is a concrete
cubic irrational number),
obtaining also exponentially small estimates
for the maximal splitting distance.
We refer to \cite{DelshamsGS04,DelshamsGG14a} for a more complete background
and references concerning exponentially small splitting, and its relation
to the arithmetic properties of the frequencies.

In this paper, we consider a 2-dimensional torus whose frequency
ratio in~(\ref{eq:omega_eps}) is given by the silver number.
Our main objective is to develop a methodology,
taking into account the arithmetic properties of the given frequencies,
allowing us to obtain asymptotic estimates for
both the maximal splitting distance
and the transversality of the splitting, as well as its continuation
for all values of $\eps\to0$.
The results on transversality and continuation generalize the results
obtained for the golden number in \cite{DelshamsG04},
and could be analogously extended to other quadratic frequency ratios
by means of a specific study in each case.

\subsection{Setup}

Here we describe the nearly-integrable
Hamiltonian system under consideration.
In particular, we study a \emph{singular} or
\emph{weakly hyperbolic} (\emph{a priori stable})
Hamiltonian with 3 degrees of freedom possessing a
2-dimensional whiskered tori with fast frequencies.
In canonical coordinates
$(x,y, \varphi, I)\in \T\times \R \times \T^2
\times \R^2$, with the symplectic form
$\df x\wedge\df y+\df\varphi\wedge\df I$, the Hamiltonian is defined by
\bea
  &&H(x,y, \varphi, I) = H_0 (x,y, I) + \mu H_1(x, \varphi),
  \label{eq:HamiltH}
\\
  &&H_0 (x, y, I) =
  \langle \omega_\eps, I\rangle + \frac{1}{2} \langle\Lambda I, I\rangle
  +\frac{y^2}{2} + \cos x -1,
  \label{eq:HamiltH0}
\\
  &&H_1 (x, \varphi)= h(x) f(\varphi).
  \label{eq:HamiltH1}
\eea
Our system has two parameters $\eps>0$ and $\mu$, linked
by a relation of kind $\mu=\eps^p$, $p>0$ (the smaller $p$ the better).
Thus, if we consider $\eps$ as the unique parameter, we have
a singular problem for $\eps\to 0$.
See \cite{DelshamsG01} for a discussion about singular and regular problems.

Recall that we are assuming a vector of fast frequencies
$\omega_\eps = \omega/\sqrt{\eps}$ as given
in~(\ref{eq:omega_eps}), with the \emph{silver frequency vector}
$\omega=(1,\Omega)$, {where $\Omega=\sqrt{2}-1$}.
It is well-known that this vector satisfies
a \emph{Diophantine condition}
\beq
|\langle k, \omega\rangle| \geq \frac{\gamma}{|k|}\,,
\;\;\; \forall k\in \Z^2\setminus\{0\}
\label{eq:DiophCond}
\eeq
with a concrete $\gamma>0$.
We also assume in~(\ref{eq:HamiltH0}) that $\Lambda$ is a symmetric
($2\times2$)-matrix, such that $H_0$ satisfies the condition of
\emph{isoenergetic nondegeneracy},
\beq
    \det \left(
    \begin{array}{cc}
    \Lambda & \omega\\
    \omega^\top & 0
    \end{array}
    \right) \neq 0.
\label{eq:isoenerg}
\eeq

For the perturbation $H_1$ in~(\ref{eq:HamiltH1}),
we consider the following periodic even functions:
\bea
  \label{eq:h}
  &&h(x) = \cos x,
\\
  \label{eq:f}
  &&f(\varphi)= \hardsum{k\in\Z^2}{k_2\ge0}
  \ee^{-\rho |k|} \cos\langle k, \varphi \rangle,
\eea
where the restriction in the sum is introduced in order to avoid repetitions.
The constant $\rho>0$
gives the complex width of analyticity of the function $f(\varphi)$.
With this perturbation, our Hamiltonian system given
by~(\ref{eq:HamiltH}--\ref{eq:f}) is \emph{reversible}
with respect to the involution
\beq\label{eq:reversible}
  \Rc:(x,y,\varphi,I)\mapsto(-x,y,-\varphi,I)
\eeq
(indeed, its associated Hamiltonian field satisfies
the identity $X_H\circ\Rc=-\Rc\,X_H$).
We point out that reversible perturbations
have also been considered in some related papers
\cite{Gallavotti94,GallavottiGM99a,RudnevW98}.
The results can be presented in a somewhat simpler way
under the assumption of reversibility.
Nevertheless, this is not essential in our approach,
and we show that our results are valid also in the
non-reversible case, if the even function $f(\varphi)$ in~(\ref{eq:f})
is replaced by a much more general function~(\ref{eq:fsigma}),
provided the phases in its Fourier expansion satisfy a suitable condition.

On the other hand, to justify the form of the perturbation $H_1$
chosen in~(\ref{eq:HamiltH1}) and~(\ref{eq:h}--\ref{eq:f}),
we stress that it makes easier the explicit
computation of the Melnikov potential, which
is necessary in order to compute explicitly the Melnikov approximation
and show that it dominates the error term in~(\ref{eq:melniapproxM}),
and therefore to establish the existence of splitting.
Moreover, the fact that all harmonics
in the Fourier expansion with respect to $\varphi$ are non-zero,
having an exponential decay, ensures
that the study of the dominant harmonics of the Melnikov potential
can be carried out directly from the arithmetic properties
of the frequency vector~$\omega$ (see Section~\ref{sect:dominant}).
It is worth reminding that the Hamiltonian defined
in~(\ref{eq:HamiltH}--\ref{eq:f}) is paradigmatic,
since it is a generalization of the famous Arnold's example
(introduced in \cite{Arnold64} to illustrate the transition
chain mechanism in Arnold diffusion).
It provides a model for the behavior of
a near-integrable Hamiltonian system near a single resonance
(see \cite{DelshamsG01} for a motivation)
and has often been considered in the literature
(see for instance \cite{GallavottiGM99b,LochakMS03,DelshamsGS04}).
Here, our aim is to emphasize the role of the arithmetic properties
of the silver frequency vector~$\omega$ in the study of the splitting.

Let us describe the invariant tori and whiskers,
as well as the splitting and Melnikov functions.
First, notice that the unperturbed system $H_0$
consists of the pendulum given by $P(x,y)= y^2/2 + \cos x -1$, and
$2$~rotors with fast frequencies:
$\dot{\varphi}= \omega_\varepsilon + \Lambda I$, $\dot{I}=0$.
The pendulum has a hyperbolic equilibrium at the origin,
and the (upper) separatrix can be parameterized by
$(x_0(s),y_0(s))=(4 \arctan\ee^s,2/\cosh s)$, $s\in\R$.
The rotors system $(\varphi, I)$ has the solutions
$\varphi = \varphi_0+(\omega_\varepsilon + \Lambda I_0)\,t$, $I= I_0$.
Consequently, $H_0$ has a $2$-parameter family of $2$-dimensional
whiskered invariant tori which have a \emph{homoclinic whisker},
i.e.~coincident stable and unstable manifolds.
Among the family of whiskered tori,
we will focus our attention on the torus located at $I=0$,
whose frequency vector is $\omega_\varepsilon$ as in~(\ref{eq:omega_eps}).

When adding the perturbation $\mu H_1$,
the \emph{hyperbolic KAM theorem} can be applied
(see for instance \cite{Niederman00})
thanks to the Diophantine condition~(\ref{eq:DiophCond})
and the isoenergetic nondegeneracy~(\ref{eq:isoenerg}).
For $\mu$ small enough, the whiskered torus persists
with some shift and deformation, as well as its local whiskers.

In general, for $\mu\ne0$ the (global) whiskers
do not coincide anymore, and one can introduce a \emph{splitting function}
giving the distance between the stable whisker $\W^\st$ and
the unstable whisker $\W^\ut$,
in the directions of the action coordinates $I\in\R^2$:
denoting $\J^{\st,\ut}(\theta)$ parameterizations
of some concrete transverse section $x=\const$ of both whiskers,
one can define the vector funcion
\ $\M(\theta):=\J^\ut(\theta)-\J^\st(\theta)$, \ $\theta\in\T^2$
\ (see \cite[\S5.2]{DelshamsG00}).
This function turns out to be the gradient
of the (scalar) \emph{splitting potential}:
\ $\M(\theta)=\nabla\Lc(\theta)$
\ (see \cite{DelshamsG00,Eliasson94}).
\ Notice that the \emph{nondegenerate critical points} of~$\Lc$ correspond
to simple zeros of $\M$ and give rise
to \emph{transverse homoclinic orbits} to the whiskered torus.

Due to the reversibility~(\ref{eq:reversible}),
the whiskers are related by the involution: $\W^\st=\Rc\,\W^\ut$.
Hence, their parameterizations can be chosen to satisfy the identity
$\J^\st(\theta)=\J^\ut(-\theta)$, provided
the transverse section $x=\pi$ is considered in their definition.
This implies that the splitting function is an odd function:
\ $\M(-\theta)=-\M(\theta)$
\ (and the splitting potential $\Lc(\theta)$ is even).
Taking into account its periodicity, we deduce that it has, at least,
the following 4~zeros (which, in principle, might be non-simple):
\beq\label{eq:4zeros}
  \theta^{(1)}_*=(0,0), \quad \theta^{(2)}_*=(\pi,0),
  \quad \theta^{(3)}_*=(0,\pi), \quad \theta^{(4)}_*=(\pi,\pi).
\eeq

Applying the Poincar\'e--Melnikov method,
the first order approximation~(\ref{eq:melniapproxM})
is given by the (vector) \emph{Melnikov function} $M(\theta)$,
which is the gradient of the \emph{Melnikov potential}:
\ $M(\theta)=\nabla L(\theta)$.
\ The latter one can be defined by integrating
the perturbation $H_1$ along
a trajectory of the unperturbed homoclinic whisker,
starting at the point of the section $s=0$ with a given phase~$\theta$:
\beq\label{eq:L}
  L(\theta) = - \int_{-\infty}^{\infty}
  [h(x_0(t))-h(0)] f(\theta +\omega_\varepsilon t)\,\df t.
\eeq
Our choice of the pendulum, whose separatrix has simple poles,
makes it possible to use the method of residues in order to compute
the coefficients of the Fourier expansion of $L(\theta)$
(see their expression in Section~\ref{sect:dominant}).
We refer to \cite{DelshamsGS04} for estimates for the Melnikov potential
and for the error term in our model~(\ref{eq:HamiltH}--\ref{eq:f}).
We stress that our approach can also be directly applied to other
classical 1-degree-of-freedom Hamiltonians $P(x,y)=y^2/2+V(x)$,
with a potential $V(x)$ having a unique nondegenerate maximum,
although the use of residues becomes more cumbersome when the separatrix
has poles of higher orders (see some examples in \cite{DelshamsS97}).

\subsection{Main result}

We show in this paper that, for the Hamiltonian
system~\mbox{(\ref{eq:HamiltH}--\ref{eq:f})}
with the 2 parameters linked by $\mu=\eps^p$,
the Poincar\'e--Melnikov method can be applied to detect the
splitting as long as we choose the exponent $p>p^*$, with some $p^*$.
Namely, we provide asymptotic estimates
for the \emph{maximal distance} of splitting, in terms of the maximum
size in modulus of the splitting function $\M(\theta)$, and for the
\emph{transversality} of the homoclinic orbits.
The main goal of this paper is to show that $\M$ has 4~simple zeros
(equivalently, that the splitting potential $\Lc$ has
4 nondegenerate critical points) for \emph{all} sufficiently small
$\eps$ and, hence, establish the
existence of 4~transverse homoclinic orbits to the whiskered tori,
generalizing the results on the \emph{continuation} of the transversality,
obtained in \cite{DelshamsG04} for the golden number.
We also obtain an asymptotic estimate
for the minimal eigenvalue (in modulus) of the splitting matrix
$\Df\M$ at each zero. This estimate provides a measure of
transversality of the homoclinic orbits.

Due to the form of $f(\varphi)$ in~(\ref{eq:f}),
the Melnikov potential $L(\theta)$
is readily represented in its Fourier series (see Section~\ref{sect:dominant}).
We use this expansion of $L$ in order to detect its
\emph{dominant harmonics} for every $\eps$.
The dominant harmonics of $L$ correspond, for $\mu$ small enough,
to the dominant harmonics of the splitting potential $\Lc$ and,
as shown in \cite{DelshamsG03},
at least 2 dominant harmonics of $\Lc$ are necessary in order
to prove the nondegeneracy of its critical points.
Such dominant coefficients
are closely related to the (quasi-)resonances of the silver frequency vector
$\omega=(1,\Omega)$.
It is established in \cite{DelshamsG03}, for any quadratic frequency vector,
a classification of the integer vectors $k$ into
\emph{primary} and \emph{secondary resonances}:
the primary resonances are the ones which fit better
the Diophantine condition~(\ref{eq:DiophCond}).
In the concrete case of the silver number $\Omega=\sqrt2-1$, the
primary resonances are related to the \emph{Pell numbers}
(see for instance \cite{FalconP07,KalmanM03}),
which play the same role as the Fibonacci numbers
in the case of the golden number considered in \cite{DelshamsG04}.
With this in mind, we define the sequence of \emph{Pell vectors}
through the following recurrence:
\beq\label{eq:silv_prim}
s_0(0)=(0,1), \;\;\; s_0(1)=(-1,2),\;\;\;  s_0(n+1)=2 s_0(n)+s_0(n-1),
\quad n\geq 1.
\eeq

We show that a change in the second dominant harmonic of
the splitting potential $\Lc$ occurs when $\eps$ goes across
some critical values $\wh\eps_n$ (called \emph{transition values}).
The nondegeneracy of the critical points of $\Lc$
can be proved in the case of 2~dominant harmonics for most values of $\eps$,
for some quadratic numbers including the silver number $\Omega$
(see \cite{DelshamsG03}).
But this excludes small neighborhoods of $\wh\eps_n$,
where the second dominant harmonic coincides with some subsequent harmonics.
In the present paper, we carry out the study near the transition values
$\wh\eps_n$ assuming that the frequency ratio $\Omega$ in~(\ref{eq:omega_eps})
is the silver number. In fact, for $\eps$ close to $\wh\eps_n$,
we need to consider 4 dominant harmonics since
the second, the third and the fourth dominant harmonics
(two of them are associated to primary resonances and one is secondary)
are of the same magnitude.
We establish, for the concrete perturbation $H_1$
in~(\ref{eq:HamiltH}--\ref{eq:f}),
the nondegenericity of the critical points of the splitting
potential $\Lc$ for values $\eps\approx\wh\eps_n$ too, and this implies
the continuation of the 4 homoclinic orbits for \emph{all} $\eps\to 0$,
with no bifurcations.

We use the notation $f \sim g$
if we can bound $c_1 |g| \leq |f| \leq c_2 |g|$ with positive constant
$c_1, c_2$ not depending on $\eps$, $\mu$.

\begin{theorem}[\emph{main result}]
  \label{thm:silvmain}
Assume for the Hamiltonian~\mbox{(\ref{eq:HamiltH}--\ref{eq:f})}
that $\eps$ is small enough
and that $\mu=\eps^p$, $p>3$.
Then, for the splitting function $\M(\theta)$ we have:
\begin{itemize}
\item[\rm(a)]
$ \displaystyle
\max_{\theta\in \T^2} |\M(\theta)|
\sim \frac{\mu}{\sqrt{\eps}}
\exp \left\{- \frac{C_0 h_1 (\eps)}{\eps^{1/4}}\right\};
$
\item[\rm(b)]
it has exactly 4 zeros $\theta^{(j)}_*$ as in~(\ref{eq:4zeros}), all simple,
and the minimal eigenvalue of $\Df\M(\theta^{(j)}_*)$ at each zero satisfies
$$
m^{(j)}_* \sim \mu \eps^{1/4}
  \exp \left\{- \frac{C_0 h_2 (\eps)}{\eps^{1/4}}\right\}.
$$
\end{itemize}
The functions $h_1(\eps)$ and $h_2(\eps)$, defined in~(\ref{eq:h12}),
are $4\ln(1+\sqrt{2})$-periodic in $\ln\eps$, with
\ $\min h_1 (\eps) =1$,
\ $\max h_1(\eps)= \min h_2(\eps) = \sqrt{(1+\sqrt{2})/2}\approx 1.0987$,
\ $\max h_2(\eps)=\sqrt{2}\approx 1.4142$.
\ On the other hand, $C_0=(\pi\rho)^{1/2}$.
\end{theorem}

In fact, we show in Section~\ref{sect:silv_epsn} that this result applies to
a much more general perturbation in~(\ref{eq:f}):
\beq\label{eq:fsigma}
f(\varphi)= \hardsum{k\in\Z^2}{k_2\ge0}
e^{-\rho |k|} \cos(\langle k, \varphi \rangle - \sigma_k),
\eeq
under a suitable condition on the phases $\sigma_k\in \T$ associated
to primary vectors $k$, see~(\ref{eq:silv_prim}).
Such a condition, established in Lemma~\ref{lm:deltatau}, will be clearly
fulfilled in our concrete reversible case~(\ref{eq:f}),
given by $\sigma_k=0$ for any $k$.

We stress that the result on continuation,
given in Theorem~\ref{thm:silvmain}, requires a careful
study of the transitions in the second dominant harmonic,
when the parameter $\eps$ goes through the values $\wh\eps_n$,
where the results of \cite{DelshamsG03} do not apply.
A result on continuation was already obtained in \cite{DelshamsG04}, but
for the golden number $\Omega=(\sqrt{5}-1)/2$, showing that, in this case,
one only needs to take into account the primary resonances.
We extend this result to the case of
the silver number $\Omega=\sqrt{2}-1$
with the additional difficulty that
at the transition values, we also have to take into account the harmonics
associated to secondary resonances.
We point out that the technique used in this paper could also be applied
to any quadratic number by means of a specific study
(in each case, assuming a suitable condition on the phases $\sigma_k$
in~(\ref{eq:fsigma})).

\bremark
If the function $h(x)$ in~(\ref{eq:h}) is replaced by~$h(x)=\cos x -1$,
then the results of Theorem~\ref{thm:silvmain} are valid for $\mu=\eps^p$
with $p>2$ (instead of $p>3$).
The details of this improvement {are not given here},
since they work exactly as in \cite{DelshamsG04}.
\eremark

\section{The silver frequency vector}\label{sect:silvfreq}

We review in this section the technique developed in \cite{DelshamsG03}
(see also \cite{DelshamsGG14b})
for studying the resonances of quadratic frequency vectors $\omega$,
in~(\ref{eq:omega_eps}), in the concrete case of the silver ratio.
This ratio has the following 1-periodic continued fraction:
$$
\Omega=\sqrt{2}-1 = [2,2,2,\ldots]
  =\frac{1}{2+\dfrac{1}{2+\dfrac{1}{2+\cdots}}}\,.
$$
It is well-known that $\omega=(1,\Omega)$, as well as
any quadratic frequency vector,
satisfies a \emph{Diophantine condition} as in~(\ref{eq:DiophCond}).
With this in mind, we define the \emph{``numerators''}
\beq\label{eq:numerators}
\gamma_k := |\langle k, \omega\rangle |\cdot|k|,
\qquad
k\in\Z^2\setminus\{0\}
\eeq
(for integer vectors, we use the norm $\abs\cdot=\abs\cdot_1$,
i.e.~the sum of absolute values of the components of the vector).
Our goal is to provide a classification of the integer vectors $k$,
according to the size of $\gamma_k$,
in order to find the primary resonances
(i.e.~the integer vectors $k$ for which $\gamma_k$ is smallest and,
hence, fitting best the Diophantine condition~(\ref{eq:DiophCond})),
and study their separation with respect to the secondary resonances.

The key point in \cite{DelshamsG03} is to use a  unimodular matrix $T$
having the vector $\omega$ as an eigenvector with eigenvalue $\lambda>1$.
This is a particular case of a result by Koch \cite{Koch99}.
For quadratic numbers, the periodicity of the continued fraction  can be
used to construct $T$ (see \cite{DelshamsGG14c,DelshamsGG14b}).
Clearly, the iterations of the matrix $T$ provide approximations
to the direction of $\omega$.
Then, the associated quasi-resonances are given
by the matrix $U:=-(T^{-1})^\top$,
according to the following important equality:
\[
  \abs{\scprod{Uk}\omega}=\frac1\lambda\abs{\scprod k\omega}.
\]
For the silver number $\Omega$, the matrices are
$$
T=\left(\begin{array}{cc}
2 & 1\\
1 & 0\\
\end{array} \right),
\qquad
U = \left(\begin{array}{cc}
0 & -1\\
-1 & 2\\
\end{array} \right).
$$
The eigenvalues of $T$ are
\beq\label{eq:lambda}
\lambda:=\Omega^{-1}=\sqrt{2}+1
\eeq
and {$-\lambda^{-1}$} with the eigenvectors $\omega=(1,\Omega)$
and $(1,-\Omega^{-1})$, respectively.
The matrix $U$ has the same eigenvectors with the eigenvalues
$-\lambda^{-1}$ and $\lambda$, respectively.
In fact, for a quadratic number \emph{equivalent} to $\Omega$,
i.e.~with a non-purely periodic continued fraction
$\hat\Omega=[b_1,\ldots,b_l, 2,2,\ldots]=[b_1,\ldots,b_l,\Omega]$,
a linear change given by a unimodular matrix can be done between
$\omega=(1,\Omega)$ and $\hat\omega=(1,\hat\Omega)$ in order
to construct the corresponding matrices $\hat T$ and $\hat U$.
This implies that the results of this paper can be extended
to any other quadratic number equivalent to $\Omega$.

We recall the results of \cite{DelshamsG03}, on the classification
of quasi-resonances for any quadratic number $\Omega$.
The study can be restricted to integer vectors
$k=(k_1,k_2)\in\Z^2\setminus\pp0$ with $\abs{\scprod k\omega}<1/2$,
and we also assume that $k_2\ge1$.
Such integer vectors have the form $k^0(j)=(-\rint(j\Omega),j)$,
where $j\ge1$ is an integer number, and $\rint(a)$
denotes the closest integer to $a$.
An integer number $j$ is said to be \emph{primitive} if
$$
\frac{1}{2\lambda} < |\langle k^0(j), \omega\rangle| < \frac{1}{2}\,.
$$
Then, the integer vectors $k\in\Z^2$ with $\abs{\scprod k\omega}<1/2$
can be subdivided into \emph{resonant sequences}:
\begin{equation}
s(j,n) := U^n k^0(j),
\qquad n=0,1,2,\ldots
\label{eq:sjn}
\end{equation}
generated by initial vectors $k^0(j)$ with a given primitive $j$.
It was proved in \cite[Th.~2]{DelshamsG03} (see also \cite{DelshamsGG14a})
that, asymptotically, each resonant sequence $s(j,n)$
exhibits a geometric growth as $n\to\infty$, with ratio $\lambda$, and that
the sequence of the numerators $\gamma_{s(j,n)}$ has a limit $\gamma^*_j$.
More precisely,
\beq\label{eq:gammajast}
|s(j,n)| = K_j \lambda^n + \Ord(\lambda^{-n}),
\qquad
\gamma_{s(j,n)} = \gamma_j^* + \Ord(\lambda^{-2n}),
\eeq
where $K_j$ and $\gamma_j^*$ can be determined explicitly for each
resonant sequence, from its primitive $j$
(see explicit formulas in \cite{DelshamsG03}).
Since the lower bounds for $\gamma^*_j$, also provided in \cite{DelshamsG03},
are increasing in $j$, we can select the minimal of them,
corresponding to some $j_0$. We denote
\beq
\gamma^*:=\liminf_{|k|\to \infty} \gamma_k
= \min_{j}\gamma^*_j = \gamma^*_{j_0}>0.
\label{eq:quad_gamj0}
\eeq
The integer vectors of the sequence $s_0(n):=s(j_0, n)$
are called \emph{the primary resonances},
and integer vectors
belonging to any of the remaining resonant sequences $s(j,n)$, $j\neq j_0$,
are called \emph{secondary resonances}.
One also introduces \emph{normalized numerators} and their limits,
after dividing by $\gamma^*$:
\[
  \tl\gamma_k:=\frac{\gamma_k}{\gamma^*}\,,
  \qquad
  \tilde{\gamma}^*_{j}:=\dfrac{\gamma^*_j}{\gamma^*}\,.
\]

For the concrete case of the silver number $\Omega=\sqrt{2}-1$, we have:
\beq\label{eq:numer_silver}
  \gamma^* = \gamma^*_1=\frac12\,,
  \qquad
  \tl\gamma_k=2\gamma_k,
  \quad
  \tl\gamma^*_j=2\gamma^*_j,
\eeq
as well as the following data, which can be obtained
from the results of \cite{DelshamsG03}:
$$
\begin{array}{llll}
j_0 = 1,  & k^0(1)= (0, 1),  & \tilde{\gamma}^*_1 =1, & K_1
= \frac{1}{2} \Omega +1 \approx 1.2071;\\[3pt]
j= 3, & k^0(3) = (-1, 3), & \tilde{\gamma}^*_3 = 2, & K_3
= \frac{3}{2} \Omega + \frac{7}{2} \approx 4.1213;\\[3pt]
j= 4, & k^0(4) = (-2, 4), & \tilde{\gamma}^*_4 = 4, & K_4
= 2\Omega +5 \approx 5.8284;\\[3pt]
j \ge 6& & \tilde{\gamma}^*_j > 6.5723 &\\[3pt]
\end{array}
$$
(notice that the integer vectors $k^0(j)$ for $j=2,5,\ldots$
are not primitive, and belong to the sequence generated by some primitive).
It is not hard to see from~(\ref{eq:sjn}), applying induction
with respect to $n$, that
$$
s(j,n)=(- p(j,n-1), p(j,n)),
$$
where $p(j, n)$ is a ``generalized'' Pell sequence:
\ $p(j, n+1) = 2 p(j, n)+ p(j, n-1)$, $n\ge1$,
starting from \ $p(j, 0) = \rint(j \Omega)$ \ and \ $p(j, 1) = j$.
\ For $j = 1$, since $\rint(\Omega) = 0$, we get the (classical) Pell sequence:
\ $P_{n+1}=2P_n+P_{n-1}$, \ with \ $P_0=0$\ and \ $P_1=1$,
and the primary resonances are
\ $s_0(n)=s(1,n)=(-P_{n-1}, P_n)$,
\ as introduced in~(\ref{eq:silv_prim}).

We denote by  $s_1(n):= s(3,n)$ the sequence of secondary vectors generated by
$k^0(3)=(-1,3)$. This sequence gives the second smallest limit
$\tilde{\gamma}^*_3=2$ and, as shown in Section~\ref{sect:silv_epsn},
it plays an essential role in
the analysis of the transversality near the transition values.
Because of this, the vectors in the sequence $s_1(n)$ will be called
\emph{the main secondary resonances}.
Using induction, we can establish the following relations between
the primary and the main secondary resonances:
\beq
s_1(n) = s_0(n)+s_0(n+1), \;\;\; n\ge0.
\label{eq:silv_s0s1}
\eeq

\section{Dominant harmonics of the splitting potential}\label{sect:dominant}

From now on, we consider the 3 degrees of freedom Hamiltonian given as
in~(\ref{eq:HamiltH}--\ref{eq:h}) but,  instead of~(\ref{eq:f}),
we consider a more general perturbation~(\ref{eq:fsigma})
with given phases $\sigma_k$.
In fact, in order to guarantee the continuation
of the transverse homoclinic orbits, a quite general condition
on the phases $\sigma_k$ will have to be fulfilled
(see this condition in~(\ref{eq:condphases})).

We put our functions $f$ and $h$ defined in~(\ref{eq:fsigma}) and~(\ref{eq:h}),
respectively,
into the integral~(\ref{eq:L}) and, calculating it by residues,
we get the Fourier expansion of the Melnikov potential:
\[
  L(\theta)=\hardsum{k\in \mathbb{Z}^2\setminus\{0\}}{k_2\ge0}
  L_k \cos(\langle k, \theta\rangle -\sigma_k),
  \qquad
  L_k = \frac{2\pi |\langle k, \omega_\varepsilon\rangle|
  \,\ee^{-\rho |k|}}{\sinh |\frac{\pi}{2}
  \langle k, \omega_\varepsilon\rangle|}\,.
\]
We point out that the phases $\sigma_k$ are the same as in~(\ref{eq:fsigma}).
Using~(\ref{eq:omega_eps}) and~(\ref{eq:numerators}),
we can present the coefficients in the form
\begin{equation}
\label{eq:alphabeta}
L_k = \alpha_k\,\ee^{- \beta_k},
\qquad
\alpha_k \approx \frac{4 \pi\gamma_k}{|k|\sqrt{\varepsilon}}\,,
\quad
\beta_k =\rho |k| + \frac{\pi \gamma_k}{2 |k|\sqrt{\varepsilon}}\,,
\end{equation}
where an exponentially small term has been neglected in the denominator
of $\alpha_k$.
For any given $\eps$, the harmonics with largest coefficients $L_k(\eps)$
correspond essentially to the smallest exponents
$\beta_k(\eps)$. Thus, we have to study the dependence on $\eps$ of
such exponents.

With this aim, we introduce for any $X$, $Y$
the function
\[
  G(\eps;X,Y):=
  \frac{Y^{1/2}}2
  \pq{\p{\frac\eps X}^{1/4}+\p{\frac X\eps}^{1/4}},
\]
having its minimum at $\varepsilon=X$,
with the minimum value $G(X;X,Y)=Y^{1/2}$.
Then, the exponents $\beta_k(\eps)$ in~(\ref{eq:alphabeta}) can be
presented in the form
\beq\label{eq:gk_quad}
  \beta_k(\eps) = \frac{C_0}{\eps^{1/4}}\,g_k (\eps),
  \qquad
  g_k (\eps):=G(\eps;\eps_k,\tl\gamma_k),
  \eeq
where
\beq\label{eq:C0}
  \eps_k:=D_0\,\frac{\tl\gamma_k^{\,2}}{\abs k^4}\,,
  \qquad
  C_0:=(\pi\rho)^{1/2},
  \qquad
  D_0:=\p{\frac{\pi}{4\rho}}^2,
\eeq
and recall that the numerators $\tl\gamma_k=2\gamma_k$
were introduced in~(\ref{eq:quad_gamj0}--\ref{eq:numer_silver}).
Consequently, for all $k$ we have
 $\beta_k(\eps)\geq \dfrac{C_0\tl\gamma_k^{1/2}}{\varepsilon^{1/4}}$\,.
This provides, according to~(\ref{eq:alphabeta}),
an asymptotic estimate for the exponent of
the maximum value of the coefficient $L_k(\eps)$ of each harmonic.

For any $\eps$ fixed we have to find the dominant terms $L_k$ and
the corresponding vectors $k$. Since the coefficients $L_k$ are exponentially
small in $\eps$,
it is more convenient to work with the functions $g_k$,
whose smallest values correspond  to the largest $L_k$.
To this aim, it is useful to consider the graphs of the functions
$g_k(\eps)$, $k\in\Z^2\setminus\{0\}$,
in order to detect the minimum of them for a given value of $\eps$.

\begin{figure}[b!]
\centering
\subfigure[]{
    \includegraphics[width=0.45\textwidth]{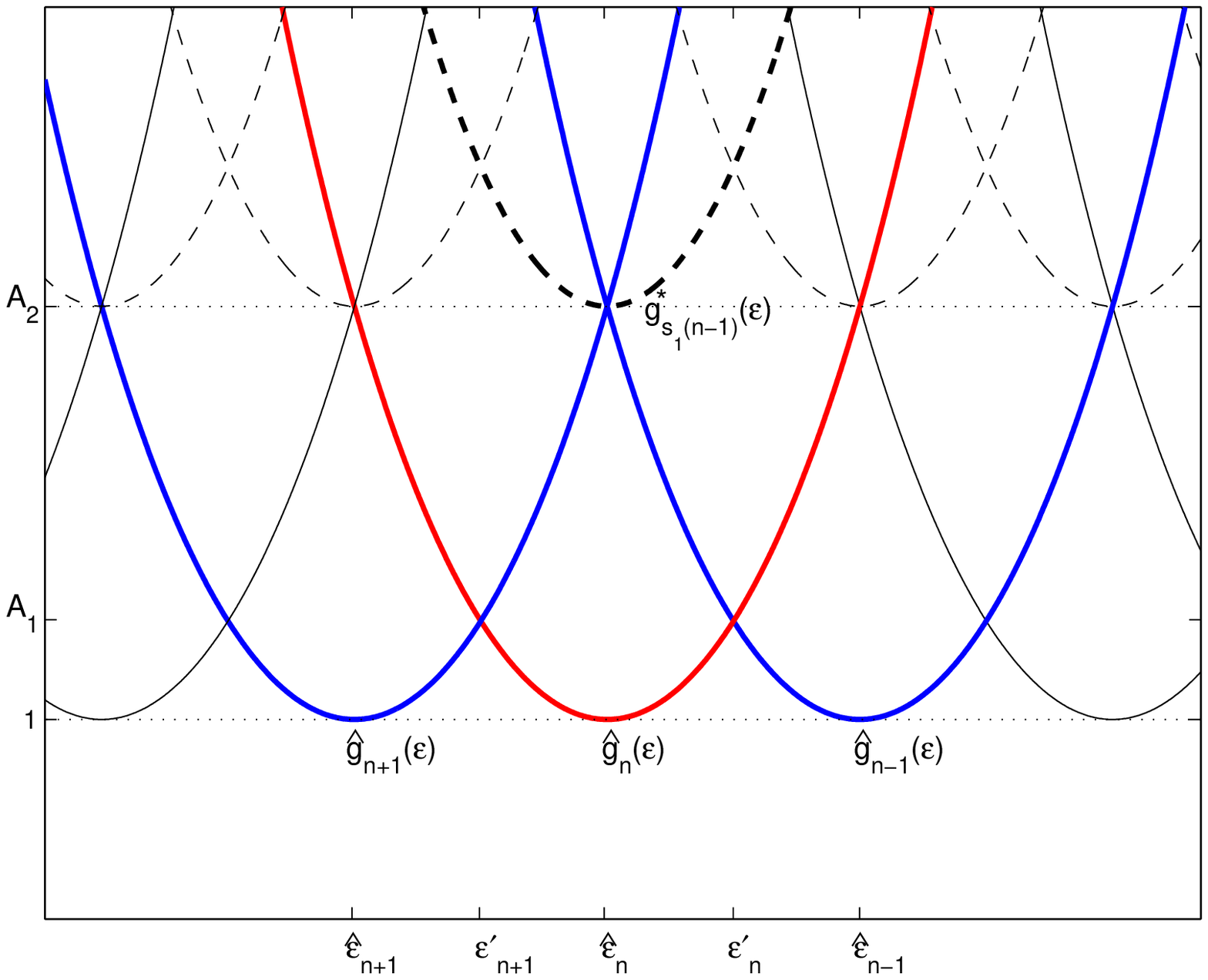}
    \label{fig:quad[2]gk}
}
\subfigure[]{
   \includegraphics[width=0.45\textwidth]{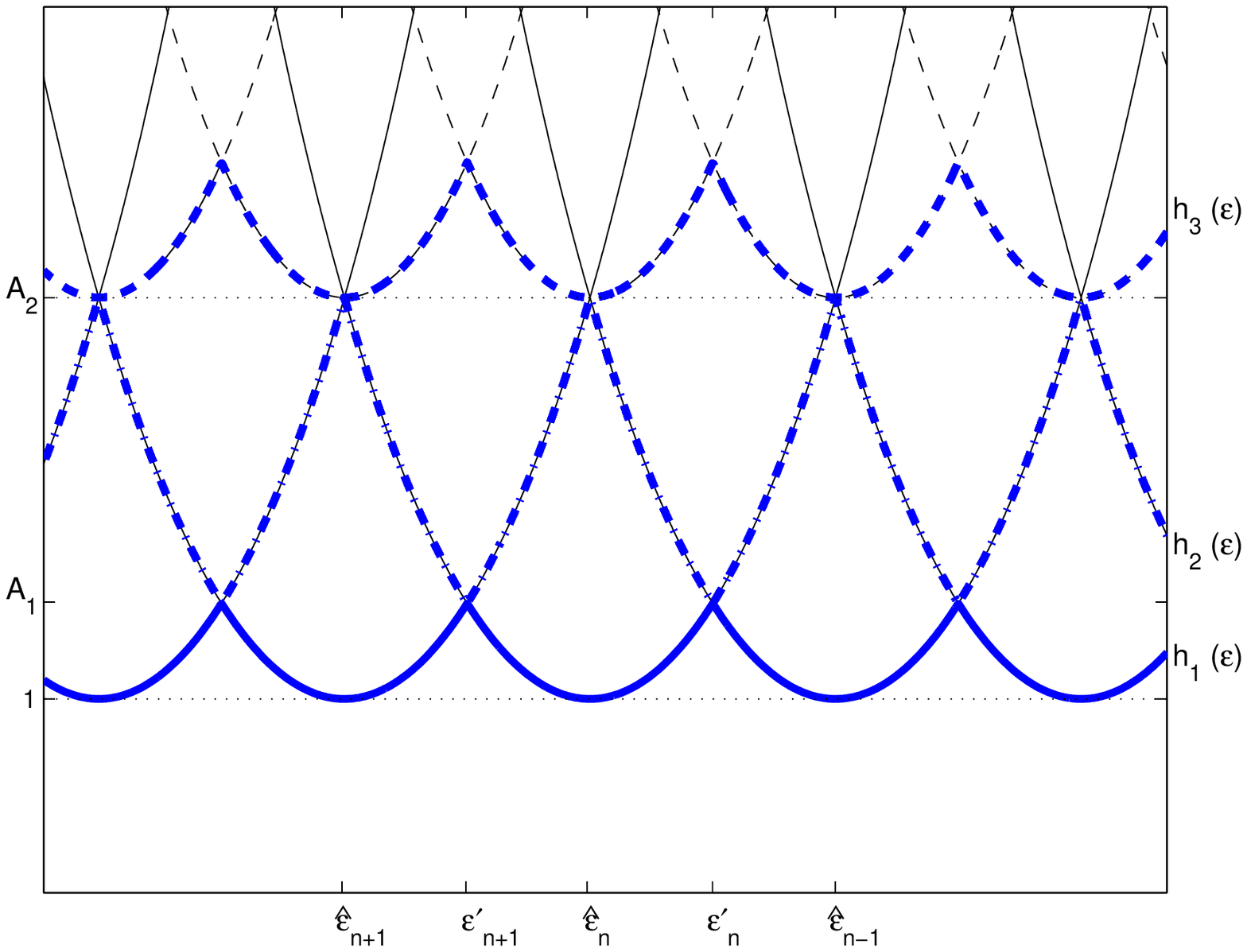}
    \label{fig:quad[2]h1}
}
\caption[]{\small
  \quad
  \subref{fig:quad[2]gk}
  \emph{Graphs of the functions $g^*_{s(j,n)}(\eps)$,
    using a logarithmic scale for $\eps$;
    the ones with solid lines are the primary functions $\wh g_n(\eps)$.}
  \quad
  \subref{fig:quad[2]h1}
  \emph{Graphs of the minimizing functions $h_1(\eps)$, $h_2(\eps)$
    and $h_3(\eps)$.}
  \quad
  \emph{Here $A_1=\sqrt{(1+\sqrt{2})/2}\approx 1.0987$
    and $A_2=\sqrt{2}\approx 1.4142$.}}
\label{fig:gkh1}
\end{figure}

We know from~(\ref{eq:gk_quad}) that the functions $g_k(\eps)$
have their minimum at $\eps=\eps_k$ and the corresponding minimal values are
$g_k(\eps_k) = \tilde{\gamma}_k^{1/2}$.
For the integer vectors $k=s(j,n)$ belonging to a resonant sequence
(recall the definition in~(\ref{eq:sjn})),
we use the approximations for $|s(j,n)|$ and
$\gamma_{s(j,n)}$ as $n\to\infty$, given in~(\ref{eq:gammajast}).
This provides the following approximations as $n\to\infty$,
$$
g_{s(j,n)}(\eps) \approx g^*_{s(j,n)}(\eps)
:=G(\eps;\eps^*_{s(j,n)},\tl\gamma^*_j),
\qquad
\eps_{s(j,n)} \approx \eps^*_{s(j,n)}
:=\frac{D_0(\tl\gamma^*_j)^2}{K_j^{\,4}\lambda^{4n}}\,.
$$
The graphs in Figure~\ref{fig:gkh1}\subref{fig:quad[2]gk},
where a logarithmic scale for $\eps$ is used,
correspond to the approximations $g^*_k(\eps)$,
rather than the true functions $g_k(\eps)$.
Note that the functions $g^*_{s(j,n)}(\eps)$  satisfy the following scaling
property:
\beq\label{eq:scaling}
  g^*_{s(j,n+1)}(\eps) = g^*_{s(j,n)} (\lambda^4 \eps).
\eeq

The case of the sequence of primary resonances
plays an important role here, since it gives the smallest minimum values
of the functions $g_k(\eps)$. With this in mind, we denote
\beq\label{eq:gn}
  \wh g_n(\eps):=g^*_{s_0(n)}=G(\eps;\wh\eps_n,1),
  \qquad
  \wh\eps_n:=\eps^*_{s_0(n)}
  =\frac{16 D_0}{\lambda^{4(n+1)}}\,,
\eeq
where we have used that $\tl\gamma^*_1=1$ and $K_1=\lambda/2$.
On the other hand, for the main secondary resonances
we can use that $\tl\gamma^*_3=2$ and $K_3/K_1=\sqrt2\,\lambda$,
and obtain
\[
  g^*_{s_1(n-1)}=G(\eps;\wh\eps_n,\sqrt2),
  \qquad
  \eps^*_{s_1(n-1)}=\wh\eps_n.
\]
Such facts are represented in Figure~\ref{fig:gkh1}\subref{fig:quad[2]gk}.

Now we define, for any given $\eps$ and for $i=1,2,3,\ldots$,
the function $h_i(\eps)$
as the $i$-th minimum of the values
$g^*_k(\eps)$, $k\in\Z^2\setminus\pp0$,
and we denote $S_i=S_i(\eps)$ the integer vectors
where {such} minima are reached:
\beq\label{eq:h12}
\begin{array}{ll}
 \displaystyle h_1(\eps):=\min_kg^*_k(\eps)=g^*_{S_1}(\eps),
  &
\displaystyle  h_2(\eps):=\min_{k\ne S_1}g^*_k(\eps) = g^*_{S_2} (\eps),
  \\
 \displaystyle h_3(\eps):=\min_{k\ne S_1, S_2}g^*_k(\eps) = g^*_{S_3} (\eps),
  &
 \displaystyle {\textrm{etc. }}
\end{array}
  \eeq
It is clear fron the scaling property~(\ref{eq:scaling}) that the
functions $h_i(\eps)$ are $4 \ln \lambda$-periodic in $\ln \eps$,
and continuous.
As we can see in Figure~\ref{fig:gkh1}\subref{fig:quad[2]h1},
the functions $h_1(\eps)$
and $h_2(\eps)$ are given by primary vectors $s_0(n)$, and $h_3(\eps)$
is given by secondary vectors $s_1(n)$.
It is easy to check that the minimum and maximum values of $h_1$ and $h_2$
are the ones given in the statement of Theorem~\ref{thm:silvmain}.

The functions $h_i(\eps)$ provide estimates, for any $\eps$, of the size
of the corresponding dominant coefficients $L_{S_i(\eps)}$
of the Melnikov potential.
We say that a given value $\eps$ is a \emph{transition value}
if $h_2(\eps)=h_3(\eps)$, since
a transition in the second dominant harmonic
takes place at these values. In the case of the silver frequencies,
these values correspond to the geometric sequence $\wh\eps_n$
defined in~(\ref{eq:gn}).
In the next section, in order to prove the tranversality
in small neighborhoods of
$\wh\eps_n$ we need to consider the 4 dominant harmonics of the splitting
potential (one of which is a main secondary resonance $s_1(n)$).
This is the main goal of this paper, since for a majority of values
of $\eps$ (excluding such neighborhoods of $\wh\eps_n$) it is enough
to consider the simpler case of 2 dominant harmonics in order to prove
the transversality, and this is already considered in \cite{DelshamsG03}
for a wider class of quadratic frequency ratios.
We also define the sequence of geometric means of the sequence $\wh\eps_n$,
\beq\label{eq:epsnprime}
\eps_n':=\sqrt{\wh\eps_n \wh\eps_{n-1}} = \frac{16 D_0}{\lambda^{4n+2}}\,,
\eeq
at which the functions $h_1(\eps)$ and $h_2(\eps)$ coincide.
For $\eps$ belonging to a given interval $(\eps'_{n+1},\eps'_n)$,
which contains the transition value $\wh\eps_n$, we have
\beq\label{eq:dominants1}
  S_1=s_0(n),
  \quad
  S_3=s_1(n+1),
\eeq
and
\beq\label{eq:dominants2}
\begin{array}{lll}
  S_2=s_0(n+1), &\quad S_4=s_0(n-1) &\quad \mbox{for \ $\eps<\wh\eps_n$,}
  \\[4pt]
  S_2=s_0(n-1), &\quad S_4=s_0(n+1) &\quad \mbox{for \ $\eps>\wh\eps_n$}
  \end{array}
\eeq
(see also Figure~\ref{fig:gkh1}).
We have the following important estimate:
since we can choose $n=n(\eps)$ such that $\eps\in(\eps'_{n+1},\eps'_n)$,
from~(\ref{eq:gammajast}) and~(\ref{eq:gn}) we obtain
\beq\label{eq:estimS}
  \abs{S_i}\sim\lambda^n\sim\eps^{-1/4},
  \quad
  i=1,2,3,4
\eeq
(recall that the notation `$\sim$' was introduced just before
Theorem~\ref{thm:silvmain}).

We will use the next lemma of \cite{DelshamsG03},
which establishes that the 4 most dominant harmonics of the Melnikov potential
are also dominant for the splitting potential,
\[
  \Lc(\theta)=\hardsum{k\in \mathbb{Z}^2\setminus\{0\}}{k_2\ge0}
  \Lc_k \cos(\langle k, \theta\rangle -\tau_k),
\]
providing an estimate for such dominant harmonics $\Lc_{S_i}$
(and an uper bound for the difference of their phases),
as well as an estimate for the sum of all other
harmonics in terms of the first neglected harmonic $\Lc_{S_5}$.
In fact, since we will be interested in some derivative of the
splitting potential, we consider the sum of (positive)
amounts of the type $|k|^l \Lc_k$.
The constant $C_0$ in the exponentials is the one defined in~(\ref{eq:C0}).

For positive amounts, we use the notation $f\preceq g$
if we can bound $f\leq c\,g$
with some constant $c$ not depending on~$\eps$ and $\mu$.

\begin{lemma}\label{lm:dominants_calL}
For $\eps$ small enough and $\mu=\varepsilon^p$ with $p>3$, one has:
\btm
\item[\rm(a)]
$\ds\Lc_{S_i}
 \sim\mu\,L_{S_i}
 \sim\frac{\mu}{\varepsilon^{1/4}}
   \,\exp\pp{-\frac{C_0h_i(\eps)}{\eps^{1/4}}}$,
\quad $\abs{\tau_{S_i}-\sigma_{S_i}}\preceq\dfrac\mu{\eps^3}$\,,
\ \ $i=1,2,3,4$;
\vspace{4pt}
\item[\rm(b)]
$\ds\sum_{k\ne S_1,\ldots, S_4}\abs{k}^l \Lc_k
 \sim\frac1{\varepsilon^{l/4}}\,\Lc_{S_5}$,
\quad $l\ge0$.
\etm
\end{lemma}

\section{Behavior near the transition values}
\label{sect:silv_epsn}

This section is devoted to the study of the transversality of the
homoclinic orbits for values of the perturbation parameter $\eps$ near the
transition values $\wh\eps_n$, defined in~(\ref{eq:gn}), where the second,
the third and the fourth dominant harmonics are of the same magnitude.
The difficulty is due to the fact that the third dominant
harmonic is associated to a main secondary resonance: $S_3=s_1(n-1)$.

We consider a concrete interval $\eps\in (\eps'_{n+1}, \eps'_n)$
which contains $\wh\eps_n$ (the values $\eps_n'$ are defined
in~(\ref{eq:epsnprime})).
For $\eps\approx\wh\eps_n$ we show that, under suitable conditions,
the splitting potential $\Lc(\theta)$ has 4 nondegenerate critical points,
which give rise to 4 transverse homoclinic orbits.
First, we study the critical points of the approximation of $\Lc(\theta)$ given
by the 4 dominant harmonics~(\ref{eq:dominants1}--\ref{eq:dominants2})
in the considered interval,
\[
  \Lc^{(4)} (\theta) :=
  \sum_{i=1,2,3,4}\Lc_{S_i}\cos(\scprod{S_i}\theta-\tau_{S_i})
\]
and, afterwards, we prove the persistence of these critical points
in the whole function $\Lc(\theta)$.

We perform the linear change
\beq
\psi_1 = \langle s_0(n-1), \theta\rangle - \tau_{s_0(n-1)},
\;\;\; \psi_2 = \langle s_0(n), \theta\rangle - \tau_{s_0(n)},
\label{eq:silv_psichange}
\eeq
that can be written as
$$
\psi = \A \theta -b, \quad \textrm{where}\ \ \A=\left(
\begin{array}{c}
s_0(n-1)^\top\\[4pt]
s_0(n)^\top
\end{array}
\right),
\quad
b= \left(
\begin{array}{c}
\tau_{s_0(n-1)}\\[2pt]
\tau_{s_0(n)}
\end{array}
\right).
$$
Since $\det\A=(-1)^{n-1}$, as easily seen from~(\ref{eq:sjn}),
this change is one-to-one on~$\T^2$.
Taking into account~(\ref{eq:silv_prim}) and~(\ref{eq:silv_s0s1}),
and recalling~(\ref{eq:dominants1}--\ref{eq:dominants2}),
we see that the function $\Lc^{(4)}(\theta)$ is transformed,
by this change, into
\bea
  \nonumber
  K^{(4)}(\psi) &= &B \cos \psi_2 + B \eta (1-Q) \cos \psi_1
  + B\eta Q\cos(\psi_1+2\psi_2 -\triangle\tau)
\\
  \label{eq:K4}
  &&+ B\eta \wtl Q\cos(\psi_1+\psi_2-\triangle\tau_1),
\eea
where we define
\bea
  \label{eq:BQQ1}
  &&B=B(\eps):=\Lc_{s_0(n)},
  \qquad
  \eta=\eta(\eps) := \displaystyle
  \frac{\Lc_{s_0(n-1)}+\Lc_{s_0(n+1)}}{\Lc_{s_0(n)}}\,.
\\
  \label{eq:BQQ2}
  &&Q=Q(\eps) := \displaystyle
  \frac{\Lc_{s_0(n+1)}}{\Lc_{s_0(n-1)}
  +\Lc_{s_0(n+1)}}\,,\quad
  \wtl Q=\wtl Q(\eps) := \displaystyle
  \frac{\Lc_{s_1(n-1)}}{\Lc_{s_0(n-1)}
  +\Lc_{s_0(n+1)}}\,,
\\
  \label{eq:deltatau}
  &&\triangle\tau:=\tau_{s_0(n+1)}-2\tau_{s_0(n)}-\tau_{s_0(n-1)},
\\
  \nonumber
  &&\triangle\tau_1:=\tau_{s_1(n-1)}-\tau_{s_0(n)}-\tau_{s_0(n-1)},
\eea
Let us describe the behavior of $Q$, $\wtl Q$, $\eta$ as $\eps$ varies
in the interval $(\eps'_{n+1}, \eps'_n)$,
which contains the transition value $\wh\eps_n$
in which we are interested.
On one hand, we see from~(\ref{eq:dominants1}--\ref{eq:dominants2})
and Lemma~\ref{lm:dominants_calL}(a)
that $\eta$ is exponentially small in $\eps$
in the whole interval, and we will consider it as a perturbation parameter.
On the other hand,
$Q$ goes from $1$ to $0$ and $\wtl Q$ takes values between $0$
and $1/2$, as $\eps$ crosses $\wh\eps_n$. More precisely,
as one can see in Figure~\ref{fig:gkh1}, for $\eps \simeq \eps_{n+1}'$
we have $\wh g_{n+1} < g_{s_1(n-1)} < \wh g_{n-1}$ and hence,
recalling~(\ref{eq:gn}),
$\Lc_{s_0(n+1)} \gg \Lc_{s_1(n-1)}\gg \Lc_{s_0(n-1)}$ and
$Q \simeq 1$, $\wtl Q \simeq 0$.
On the other hand,  for $\eps \simeq \eps_{n}'$
we have $\wh g_{n+1} > g_{s_1(n-1)} > \wh g_{n-1}$ and hence
$Q\simeq 0$, $\wtl Q \simeq 0$.
At $\eps=\wh\eps_n$ we have $\wh g_{n+1}=g_{s_1(n-1)}=\wh g_{n-1}$,
and therefore the harmonics coincide and we have $Q=\wtl Q = 1/2$.
In the interval $(\eps'_{n+1}, \eps'_n)$ considered, we see that
$Q$ is decreasing, and $\wtl Q$ has a maximum at $\wh\eps_n$ and lies
between $Q$ and $1-Q$.

We are going to use the following lemma, whose proof is a simple application
of the standard fixed point theorem.

\begin{lemma}
If $F:\T\to\R$ is differentiable and satisfies
\ $(F')^2 + F^2 < 1$,
\ then the equation
\ $\sin x = F(x)$
\ has exactly two solutions $\ol x$ and $\ool x$, which are simple.
Furthermore, if $F(x)=\Ord(\eta)$ for any $x\in \T$ with $\eta$ sufficiently
small, then the solutions of the equation satisfy
$\ol x = \Ord(\eta)$ and $\ool x=\pi+\Ord(\eta)$.
\label{lm:FPsin1D}
\end{lemma}

Now, we introduce the following important quantity:
\bea
  \label{eq:silv_Eastr}
  &&E^*= E^*(\eps) := \min (E^{(+)}, E^{(-)}),
  \quad\textrm{where}
\\
  \nonumber
  &&E^{(\pm)} := \sqrt{\left[1-Q + Q \cos \triangle\tau \pm \wtl Q
\cos \triangle\tau_1\right]^2 + \left[Q \sin \triangle\tau \pm \wtl Q
\sin \triangle\tau_1\right]^2}\,.
\eea
In the next lemma we prove the existence of 4 critical points
of $K^{(4)}$ for $\eta$ small enough, provided $E^*>0$.

\begin{lemma} Assume that, in~(\ref{eq:silv_Eastr}),
\beq\label{eq:EasstrPosit}
E^*(\eps)>0, \qquad \forall \eps \in (\eps'_{n+1}, \eps'_{n}).
\eeq
If  $\eta \preceq E^*$ in~(\ref{eq:BQQ1}),
the function $K^{(4)}(\psi)$ introduced in~(\ref{eq:K4}) has
4~nondegenerate critical points
\ $\psi^{(j)}=\psi^{(j),0}+\Ord(\eta)$, $j=1,2,3,4$,
where we define
\beq\label{eq:4points}
  \begin{array}{ll}
    \psi^{(1),0}=(\alpha^{(+)},0),
    &\quad\psi^{(2),0}=(\alpha^{(+)}+\pi,0),
  \\[4pt]
    \psi^{(3),0}=(\alpha^{(-)},\pi),
    &\quad\psi^{(4),0}=(\alpha^{(-)}+\pi,\pi),
  \end{array}
\eeq
with
\beq
\cos \alpha^{(\pm)}
= \frac{1-Q + Q \cos \triangle\tau \pm \wtl Q \cos \triangle\tau_1}
{E^{(\pm)}}\,,
\quad
\sin \alpha^{(\pm)} = \frac{Q \sin \triangle\tau \pm \wtl Q
\sin \triangle\tau_1}{E^{(\pm)}}\,.
\label{eq:silv_alppm}
\eeq
At the critical points,
\bean
&&|\det \Df^2K^{(4)} (\psi^{(1,2)})|= B^2 \eta (E^{(+)} + \Ord(\eta)),
\\
&&|\det \Df^2K^{(4)} (\psi^{(3,4)})|= B^2 \eta (E^{(-)} + \Ord(\eta)).
\eean
\label{lm:silvK4}
\end{lemma}

\proof
The critical points of $K^{(4)}(\psi)$ are the solutions
to the system of equations
\beq
\begin{array}{l}
(1-Q) \sin \psi_1 +  Q \sin (\psi_1 + 2 \psi_2 -\triangle\tau)
+ \wtl Q \sin (\psi_1 + \psi_2 - \triangle\tau_1)=0,\\[4pt]
\sin \psi_2 + 2  \eta Q \sin (\psi_1 + 2 \psi_2 - \triangle\tau)
+ \eta \wtl Q \sin (\psi_1 + \psi_2 - \triangle\tau_1) =0.
\end{array}
\label{eq:silv_diffK4}
\eeq
We can rewrite the second equation as follows:
\bea
  \label{eq:silv_f}
  &&\sin \psi_2 = \eta f(\psi_1, \psi_2),
\\
  \nonumber
  &&\textrm{where}\quad
  f(\psi_1,\psi_2) := -2 Q \sin (\psi_1 + 2 \psi_2 - \triangle\tau)
  -\wtl Q \sin (\psi_1 + \psi_2 - \triangle\tau_1).
\eea
Since $\eta$ is small enough and $f$ is bounded with its derivatives,
we can apply Lemma~\ref{lm:FPsin1D} with $F=\eta f$,
and $\psi_1$ as a parameter,
and we get that equation~(\ref{eq:silv_f}) has two solutions:
\ $\ol\psi_2= \ol\psi_2 (\psi_1) = \Ord(\eta)$ \ and
\ $\ool\psi_2= \ool\psi_2 (\psi_1)= \pi + \Ord(\eta)$.

Substituting $\ol\psi_2(\psi_1)$
into the first equation of~(\ref{eq:silv_diffK4}),
we get an equation $F^{(+)}_\eta (\psi_1) = 0$, with the function
\bean
F^{(+)}_\eta  &:= &(1-Q) \sin \psi_1 + Q \sin (\psi_1-\triangle\tau)
+ \wtl Q\sin(\psi_1-\triangle\tau_1)\\
&&-\eta f^{(+)}(\psi_1,\ol\psi_2;\eta)\\
&= &\left[1-Q + Q \cos \triangle\tau
+ \wtl Q \cos \triangle\tau_1\right] \sin \psi_1\\
&&- \left[Q \sin \triangle\tau + \wtl Q \sin \triangle\tau_1\right]
\cos \psi_1- \eta f^{(+)} (\psi_1,\ol\psi_2; \eta)\\
&= &E^{(+)} \sin (\psi_1 - \alpha^{(+)})-\eta f^{(+)}(\psi_1,\ol\psi_2;\eta),
\eean
where $E^{(+)}$ and $\alpha^{(+)}$ are the constants
defined in~(\ref{eq:silv_Eastr}) and~(\ref{eq:silv_alppm}),
respectively, and a function
$f^{(+)}$, which is bounded jointly with its derivatives.
Thus, provided $E^{(+)} > 0$,
the equation $F^{(+)}_\eta = 0$ is {equivalent} to
$$
\sin (\psi_1 - \alpha^{(+)} )=
\frac{\eta}{E^{(+)}}\,f^{(+)}(\psi_1,\ol\psi_2;\eta)
$$
and, by Lemma~\ref{lm:FPsin1D}, it
has 2 solutions $\psi^{(1)}_1 = \alpha^{(+)} +  \Ord(\eta)$
and $\psi^{(2)}_1 = \alpha^{(+)} + \pi + \Ord(\eta)$,
since $\eta \preceq E^*\le E^{(+)}$.
In this way, we have 2 critical points as solutions of
the system~(\ref{eq:silv_diffK4}):
\ $\psi^{(j)}=(\psi^{(j)}_1,\ol\psi_2(\psi^{(j)}_1))$, $j=1,2$.

We proceed analogously for $\ool\psi_2$ and rewrite the first equation
of~(\ref{eq:silv_diffK4}) as
$$
F^{(-)}_\eta :=
E^{(-)} \sin (\psi_1 - \alpha^{(-)})-\eta f^{(-)}(\psi_1,\ool\psi_2;\eta)=0.
$$
Assuming that $E^{(-)} > 0$,
we get other two solutions $\psi^{(3)}_1 = \alpha^{(-)} + \Ord(\eta)$
and $\psi^{(4)}_1= \alpha^{(-)} +\pi + \Ord(\eta)$,
since $\eta \preceq E^*\le E^{(-)}$.
Such solutions give rise to the other 2 critical points
$\psi^{(j)}$, $j=3,4$.

To compute the determinant at the critical points, we use that
\bean
\lefteqn{\det \Df^2K^{(4)} (\psi)}
\\
&= &B^2(\eta[(1-Q) \cos\psi_1
+  Q \cos(\psi_1 + 2 \psi_2 -\triangle\tau)
\\
&&+ \wtl Q\cos (\psi_1 + \psi_2 - \triangle\tau_1)]
\cdot\cos \psi_2
+ \Ord(\eta^2))
\eean
for any $\psi\in \T^2$. At $\psi^{(1)}$, for example, we have
\bean
\det \Df^2K^{(4)} (\psi^{(1)})
&= &B^2 \left(\eta \left.\frac{\partial F^{(+)}_\eta}
{\partial \psi_1}\right|_{\psi^{(1)}}\cdot\cos \psi^{(1)}_2
+ \Ord(\eta^2)\right)
\\
&=&B^2 (\eta E^{(+)} + \Ord(\eta^2)),
\eean
and similarly with the other 3 critical points.
\qed

\bremark
In our case of a reversible perturbation, as introduced in~(\ref{eq:f}),
we obtain in~(\ref{eq:silv_alppm}) the values $\alpha^{(\pm)}=0$.
By the linear change~(\ref{eq:silv_psichange}),
and using that the phases are $\sigma_k=0$,
we get the 4~critical points of $\Lc^{(4)}$,
as deduced in~(\ref{eq:4zeros}) from the reversibility property.
\eremark

To ensure the existence of nondegenerate critical points of $K^{(4)}$,
in Lemma~\ref{lm:silvK4} we have assumed condition~(\ref{eq:EasstrPosit}).
In the next lemma we see when this assumption fails.

\begin{lemma}\label{lm:Eastr}
Let $0<Q<1$, $0<\wtl Q \leq 1/2$ and
$\triangle\tau,\triangle\tau_1\in \mathbb{T}$ given,
and consider $E^*$ defined as in~(\ref{eq:silv_Eastr}).
Then, one has $E^*=0$ if and only if
the following three conditions are satisfied:
\begin{equation}\label{eq:costau}
 |1-2Q|\leq \wtl Q, \quad \cos  \triangle\tau
=  -\frac{1-2 Q + 2 Q^2 -\wtl Q^2}{2(1-Q) Q}\,,\quad
\cos \triangle\tau_1=\pm\frac{\wtl Q^2+1-2 Q}{2(1-Q) \wtl Q}\,.
\end{equation}
\end{lemma}

\begin{figure}[b!]
\centering
\subfigure[]{
    \includegraphics[width=0.5\textwidth]{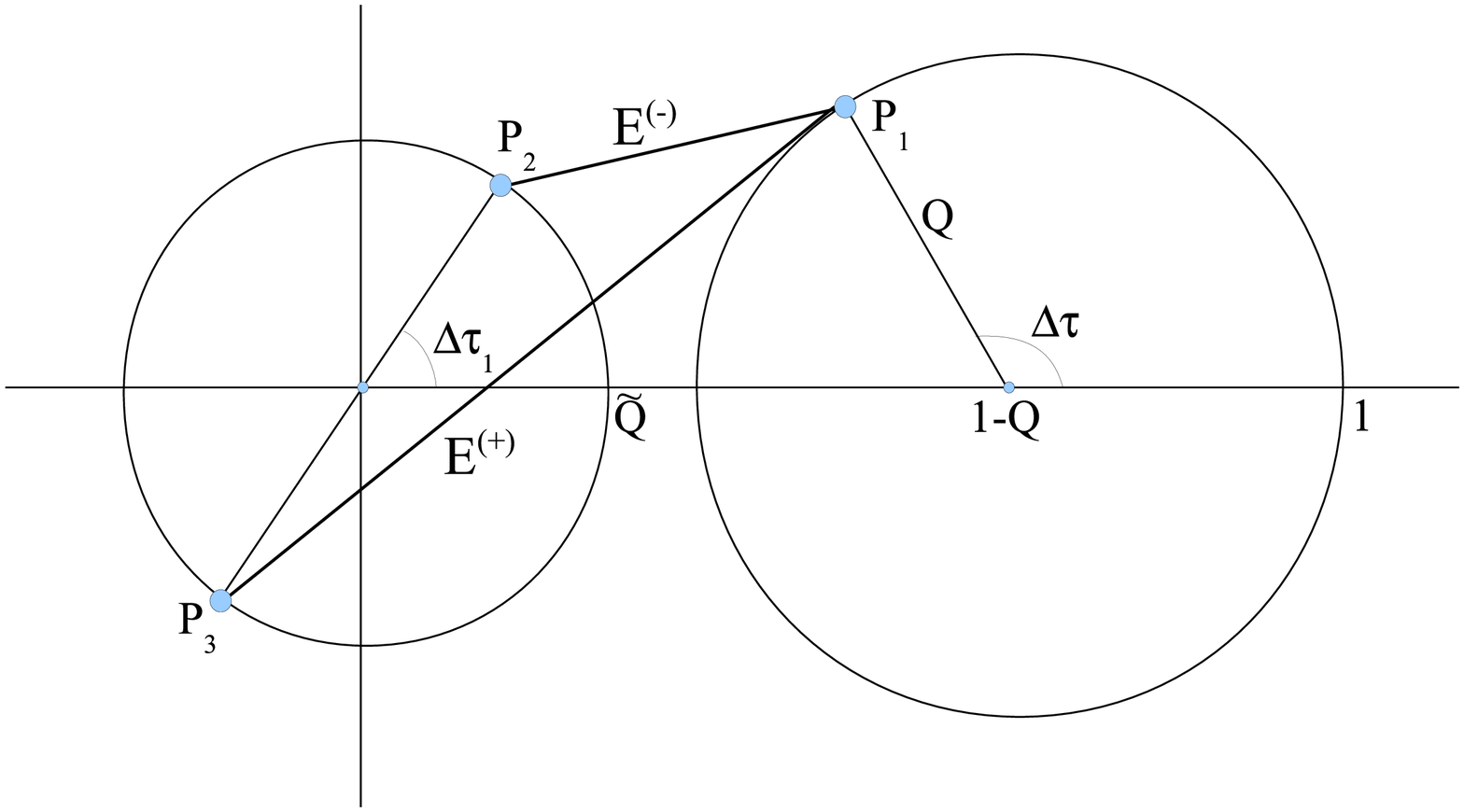}
    \label{fig:Epm}
}\qquad
\subfigure[]{
   \includegraphics[width=0.4\textwidth]{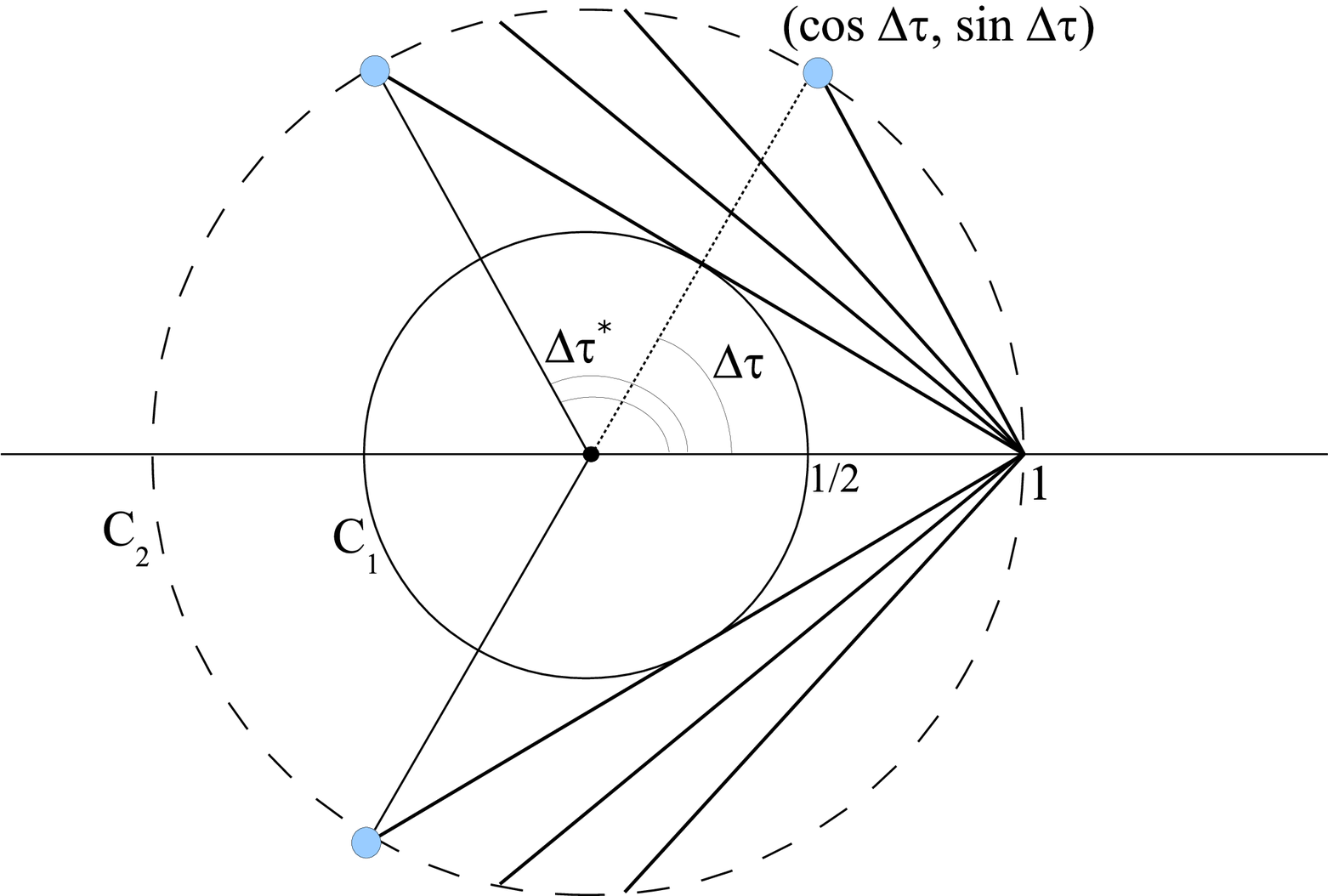}
   \label{fig:Epmneq0}
}
\caption[]{\small
  \quad
  \subref{fig:Epm}
  \emph{Geometrical representation of $E^{(+)}$ and $E^{(-)}$.}
  \\
  \subref{fig:Epmneq0}
  \emph{$E^{(\pm)}>0$ (the straight lines do not intersect
    the circle $C_1$) if $|\triangle\tau|<\triangle\tau^*=2\pi/3$.}}
\label{fig:EpmGeom}
\end{figure}

\proof We prove this lemma geometrically.
It is clear from~(\ref{eq:silv_Eastr})
that $E^*=0$ if and only if $E^{(+)}=0$ or $E^{(-)}=0$,
i.e.~one of the following two assertions hold:
$$
\begin{array}{lll}
1-Q + Q \cos \triangle\tau
= - \wtl Q \cos \triangle\tau_1 & \textrm{ and } & Q \sin \triangle\tau
=- \wtl Q \sin \triangle\tau_1, \\[4pt]
1-Q + Q \cos \triangle\tau
= \wtl Q \cos \triangle\tau_1 & \textrm{ and } & Q \sin \triangle\tau
= \wtl Q \sin \triangle\tau_1.
\end{array}
$$
Now, we consider the points
\bean
  &P_1= (1-Q +Q \cos \triangle\tau, Q \sin \triangle\tau),
\\[4pt]
  &P_2 = (\wtl Q \cos \triangle\tau_1, \wtl Q \sin \triangle\tau_1),
  \quad
  P_3 = (-\wtl Q \cos \triangle\tau_1, -\wtl Q \sin \triangle\tau_1),
\eean
which lie on the circles represented in
Figure~\ref{fig:EpmGeom}\subref{fig:Epm}, and, hence,
$E^{(+)}$ is the distance $P_1 P_3$, while $E^{(-)}$ is $P_1 P_2$.

Varying $Q$ and $\wtl Q$ and changing the corresponding circles
in Figure~\ref{fig:EpmGeom}\subref{fig:Epm} in order to see
when $P_1$ coincides
with $P_2$ and, thus, $E^{(-)}=0$,
one can get that if $|1-2Q|<\wtl Q$, the circles intersect
(at the point $P_1\equiv P_2$)
and there is a triangle
with sides $Q$, $1-Q$, $\wtl Q$ and angles satisfying~(\ref{eq:costau}).
For $1-2Q=\pm \wtl Q$, the circles are tangent having $\triangle\tau =\pi$,
and $\triangle\tau_1=0$ (if $1-2Q=\wtl Q$)
or $\triangle\tau_1=\pi$ (if $1-2Q=-\wtl Q$).
In the case $|1-2Q|>\wtl Q$, the circles do not intersect.
The case $E^{(+)}=0$ (which corresponds to $P_1\equiv P_3$)
can be studied in a similar way.
\qed

In this way, we can ensure the existence and continuation
of the 4~critical points of $K^{(4)}$
given by Lemma~\ref{lm:silvK4} if the three conditions~(\ref{eq:costau})
do no hold simultaneously for any $\eps\in(\eps'_{n+1},\eps'_n)$.
Now, we provide a simple \emph{sufficient condition} allowing us to avoid
the occurrence of~(\ref{eq:costau}) and, hence,
to ensure~(\ref{eq:EasstrPosit}).

\begin{lemma}\label{lm:deltatau}
If
\beq
|\triangle\tau| < \frac{2\pi}{3}\,,
\label{eq:Epos}
\eeq
then the condition~(\ref{eq:EasstrPosit}) is fulfilled
independently of $Q$, $\wtl Q$, $\triangle\tau_1$.
\end{lemma}

\proof
This is a corollary of Lemma~\ref{lm:Eastr}.
Indeed, the inequality~(\ref{eq:Epos}) implies that
we have $\cos \triangle\tau > -1/2$.
Then, if the second equality in~(\ref{eq:costau}) is satisfied,
we have $1-\wtl Q^2 < 3Q(1-Q)$,
which contradicts the facts that
$0<Q<1$ and $0<\wtl Q \leq 1/2$.
\qed

\bremark We can provide a geometric interpretation for this lemma.
In Figure~\ref{fig:EpmGeom}\subref{fig:Epmneq0} we consider two circles
centered at the origin: $C_1$ with radius
$1/2$ (the maximum value for $\wtl Q$) and the unit circle $C_2$.
For any given $\triangle\tau$,
the map $Q \mapsto (1- Q+Q \cos \triangle\tau, Q \sin \triangle\tau)$,
for $0\leq Q\leq 1$, gives us a family of straight lines
(with $\triangle\tau$ as a parameter) connecting the
points $(1,0)$ and $(\cos \triangle\tau, \sin \triangle\tau)$,
both belonging to $C_2$.
The straight lines corresponding to $\triangle\tau$ satisfying~(\ref{eq:Epos})
do not intersect the circle $C_1$, which implies that $E^*> 0$
(see the proof of Lemma~\ref{lm:Eastr}).
Notice that, for $\wtl Q=1/2$, the critical value $\triangle\tau^*=2\pi/3$
is sharp, but for $\wtl Q< 1/2$ the critical value would be
greater: $\triangle\tau^*>2\pi/3$.
\eremark

In the next lemma, we prove the persistence of the 4~critical points
$\psi^{(j)}$ of the approximation $K^{(4)}(\psi)$,
when the non-dominant terms are also considered.
With this aim, we denote $K(\psi)$ the function obtained
when the linear change~(\ref{eq:silv_psichange})
is applied to the \emph{whole} splitting potential $\Lc(\theta)$.
Recalling the definitions~(\ref{eq:BQQ1}--\ref{eq:BQQ2}), we can write:
$$
K(\psi) = K^{(4)}(\psi) + B \eta \eta' {G}(\psi),
$$
where the term $B \eta \eta' {G}(\psi)$ corresponds to the sum of all
non-dominant harmonics, and $\Lc_{S_5} = B \eta \eta'$ is the
largest among them with
\beq\label{eq:etaprime}
\eta':=\frac{\Lc_{S_5}}{\Lc_{s_0(n-1)}
+\Lc_{s_0(n+1)}} \ll Q, \wtl Q.
\eeq
Note that the function $G$ is obtained via the
linear change~(\ref{eq:silv_psichange}) applied to the non-dominant
harmonics of $\Lc(\theta)$. Thus, using
Lemma~\ref{lm:dominants_calL}(b), we get bounds
for $G(\psi)$ and its partial derivatives:
\beq
|G|\preceq 1, \quad
|\partial_{\psi_i} G| \preceq \eps^{-1/2}, \quad
|\partial^{\,2}_{\psi_i \psi_j} G| \preceq \eps^{-1}, \quad i,j =1,2,
\label{eq:silv_bndG}
\eeq
where we have taken into account that, by~(\ref{eq:estimS}),
the entries of the matrix of the linear change~(\ref{eq:silv_psichange})
are $\sim\eps^{-1/4}$.

\begin{lemma}
Assuming condition~(\ref{eq:EasstrPosit}),
if $\bar\eta:=\max(\eta,\eta\eta'\eps^{-1})\preceq E^*$, then the
function $K(\psi)$ has 4 critical points, all nondegenerate:
 $\psi_*^{(j)}= \psi^{(j),0} + \Ord(\bar\eta)$, $j=1,2,3,4$,
with $\psi^{(j),0}$ defined in~(\ref{eq:4points}).
At the critical points,
\bean
&&|\det \Df^2K (\psi_*^{(1,2)})|= B^2 \eta (E^{(+)}+ \Ord(\bar\eta)),
\\
&&|\det \Df^2K (\psi_*^{(3,4)})|= B^2 \eta (E^{(-)}+ \Ord(\bar\eta)).
\eean
\label{lm:silvKtot}
\end{lemma}
\proof
The critical points of $K(\psi)$ are the solution of the following equations,
which are perturbations of~(\ref{eq:silv_diffK4}):
\bean
&&(1-Q) \sin \psi_1 +  Q \sin (\psi_1 + 2 \psi_2 -\triangle\tau)
+ \wtl Q \sin (\psi_1 + \psi_2 - \triangle\tau_1)\\
&&\hspace{9cm}
- \eta \eta' \partial_{\psi_1}G=0,\\
&&\sin \psi_2 + 2  \eta Q \sin (\psi_1 + 2 \psi_2 - \triangle\tau)
+\eta \wtl Q \sin (\psi_1 + \psi_2 - \triangle\tau_1)
- \eta \eta' \partial_{\psi_2}G  =0.
\eean
Now, we can proceed as in the proof of Lemma~\ref{lm:silvK4}.
Indeed, applying Lemma~\ref{lm:FPsin1D} twice
we can solve the second equation for $\psi_2$ with $\psi_1$ as a parameter,
and we replace the solution in the first equation and solve it for $\psi_1$.
The only difference with respect to Lemma~\ref{lm:silvK4}
is that now we have additional perturbative terms
$\eta\eta'\partial_{\psi_i}G$, which we have bounded in~(\ref{eq:silv_bndG}),
and for this reason we consider $\bar\eta$ as the size of the perturbation.
The determinant at the critical points can be computed
as in Lemma~\ref{lm:silvK4}.
\qed

\bremark
The smallness condition on $\bar\eta$ in Lemma~\ref{lm:silvKtot}
is clearly fulfilled in our case, since in~(\ref{eq:etaprime})
we have that $\eta'$ is exponentially small in $\eps$ and,
hence, can be bounded by any power of $\eps$.
\eremark

Applying the inverse (one-to-one) of
the linear change~(\ref{eq:silv_psichange}), the 4~critical points
$\psi_*^{(j)}$ of $K(\psi)$ give rise to 4~critical points of $\Lc(\theta)$,
all nondegenerate:
\beq
\theta^{(j)}_* = \A^{-1} (\psi_*^{(j)} + b), \;\;\; j=1, 2, 3, 4.
\label{eq:silv_crptheta}
\eeq

\begin{lemma}
Assuming condition~(\ref{eq:EasstrPosit}),
if $\bar\eta:=\max(\eta,\eta\eta'\eps^{-1})\preceq E^*$, then
the splitting potential $\Lc(\theta)$ has exactly 4  critical
points $\theta^{(j)}_*$, given by~(\ref{eq:silv_crptheta}), all nondegenerate,
and the minimal eigenvalue (in modulus) $m^{(j)}_*$ of
$\Df^2\Lc(\theta^{(j)}_*)$ satisfies
\[
E^* \sqrt{\eps}\,\Lc_{S_2} \preceq m^{(j)}_*
\preceq\sqrt{\eps}\,\Lc_{S_2}, \;\;\; j=1, 2, 3, 4.
\]
\label{lm:silv_critpL}
\end{lemma}

\proof
The proof is similar to the one of \cite[Lemma 5]{DelshamsG03} and, thus,
we give here only a sketch of the proof.
First, denoting $D = \det \Df^2\Lc(\theta^{(j)}_*)$ and
$T = \tr \Df^2\Lc(\theta^{(j)}_*)$, it is not hard to see that,
if $\abs D\ll T^2$, then $m^{(j)}_*\sim |D|/|T|$.
Thus, we need to provide asymptotic estimates for $\abs D$ and $\abs T$.

Since $|\det \A| = 1$, the matrices $\Df^2K(\psi_*^{(j)})$ and
$\Df^2\Lc(\theta^{(j)}_*)
= \A^{\top} \Df^2K(\psi_*^{(j)}) \A$ have equal determinants,
and, hence, by Lemma~\ref{lm:silvKtot},
\[
|D| = B^2 \eta (E^{(\pm)} + \Ord(\bar\eta))
\sim E^{(\pm)} \Lc_{s_0(n)} (\Lc_{s_0(n-1)}
+\Lc_{s_0(n+1)})
\sin
E^{(\pm)} \Lc_{S_1} \Lc_{S_2},
\]
where we have taken into account the definitions~(\ref{eq:BQQ1})
and the relations~(\ref{eq:dominants1}--\ref{eq:dominants2})
between the dominant harmonics and the primary resonances.
Using that $E^*\le E^{(\pm)}\preceq 1$, we get a lower and an upper bound
for $\abs D$.

On the other hand, for the components of
$\Df^2K (\psi_*^{(j)}) = \left( \begin{array}{cc} k_{11} & k_{12}\\ k_{12}
& k_{22}\end{array} \right)$,
given in first approximation by derivatives of~(\ref{eq:K4}),
we have $\abs{k_{22}}\sim B(1+\Ord(\bar\eta))$ as the main entry,
and $\abs{k_{11}},\abs{k_{12}}\preceq B\bar\eta$.
By the linear change~(\ref{eq:silv_psichange})
the trace of $\Df^2\Lc(\theta^{(j)}_*)$ is given by
$$
T = k_{11} \langle s_0(n-1), s_0(n-1)\rangle
+ 2 k_{12} \langle s_0(n-1),s_0(n)\rangle
+ k_{22} \langle s_0(n), s_0(n) \rangle.
$$
Then, applying~(\ref{eq:estimS}) and the estimates
of Lemma~\ref{lm:dominants_calL}(a), we obtain
$$
|T| \sim \frac{1}{\sqrt{\eps}}\,\Lc_{S_1}.
$$
Now, we have an estimate for the quotient $\abs D/\abs T$,
which gives us the desired estimate for the minimal eigenvalue.
\qed

\proofof{Theorem~\ref{thm:silvmain}}
Finally, we can complete the proof of our main result.
As explained in Section~\ref{sect:dominant}, to establish the transversality
for all sufficiently small $\eps$, it is enough to consider a neighborhood
of the transition values $\wh\eps_n$, since for other values of $\eps$
it is enough to consider 2~dominant harmonics and the results
of \cite{DelshamsG03} apply.

For $\eps$ close to a transition value $\wh\eps_n$,
recalling that $\M(\theta)
= \nabla \Lc(\theta)$, it follows from
Lemma~\ref{lm:silv_critpL} that, under~(\ref{eq:EasstrPosit}),
the splitting function
$\M(\theta)$ has 4 simple zeros $\theta_*$, given in~(\ref{eq:silv_crptheta}).
Likewise, by Lemma~\ref{lm:deltatau} the condition~(\ref{eq:EasstrPosit})
is fulfilled if
\beq\label{eq:condphases}
|\sigma_{s_0(n+1)}-2\sigma_{s_0(n)}-\sigma_{s_0(n-1)}|\approx|\triangle\tau|
< \frac{2 \pi}{3}\,,
\qquad
\forall n\ge1
\eeq
(we have taken into account the bound on the difference of phases
$\sigma_k$ and $\tau_k$ given in Lemma~\ref{lm:dominants_calL}(a)).
The particular case of a reversible perturbation~(\ref{eq:f})
corresponds to~(\ref{eq:fsigma}) with $\sigma_k =0$ for every $k$,
and hence condition~(\ref{eq:condphases}) on the phases is clearly fulfilled.
Moreover, we have $E^{(\pm)}=1\pm \wtl Q\sim1$ in~(\ref{eq:silv_Eastr}),
and hence $E^*=1-\wtl Q\ge1/2$, which implies that $1/2\le E^*\le1$.
By Lemma~\ref{lm:silv_critpL}, for the minimal eigenvalue
of the splitting matrix $\Df\M(\theta_*)$ at each zero we can write
$m_*\sim\sqrt\eps\,\Lc_{S_2}$.
This estimate, together with the estimate
on $\Lc_{S_2}$ given by Lemma~\ref{lm:dominants_calL}, implies part~(b).

As for part~(a),
the maximal splitting distance is given by the most dominant harmonic
$$\max\limits_{\theta\in \mathbb{T}^2} |\M(\theta)| \sim \abs{\M_{S_1}}
 \sim\mu\abs{S_1}L_{S_1}$$
(see for instance \cite{DelshamsGG14a}),
and the corresponding estimate of Lemma~\ref{lm:dominants_calL} implies the
desired estimate.
\qed

\bremark
For the sake of simplicity, we have restricted the statement of 
Theorem~\ref{thm:silvmain} to the case of a reversible perturbation
given by~(\ref{eq:f}) with the phases $\sigma_k=0$.
Nevertheless, our results apply to a much more
general perturbation~(\ref{eq:fsigma}),
provided the phases $\sigma_{s_0(n)}$,
associated to the primary resonances,
satisfy the inequality~(\ref{eq:condphases}).
\eremark

\small

\def\noopsort#1{}

\end{document}